\documentclass{article}

\input{tcilatex}

\begin{document}

\title{A Geometric Construction for the Evaluation of Mean Curvature}
\author{Pavel Grinfeld}
\date{7/20/2011}
\maketitle

\begin{abstract}
We give a relationship that yields an effective geometric way of evaluating
mean curvature of surfaces. The approach is reminiscent of the Gauss's
contour based evaluation of intrinsic curvature. The presented formula may
have a number of potential applications including estimating the normal
vector and mean curvature on triangulated surfaces.

Given how brief is its derivation, it is truly surprising that this formula
does not appear in the existing literature on differential geometry -- at
least according to the author's search. We hope to learn about a reference
containing this result.
\end{abstract}

The presented formula admits a number of straightforward generalizations
including to higher dimensions and submanifolds of co-dimension other than
one. Here, we discuss the formula for a two dimensional surface embedded in
a three dimensional Euclidean space.

Consider a regular patch $P$ within a smooth surface $S$. Suppose that the
area of $P$ is $A$ and that its smooth contour boundary is $\Gamma $. Let $%
\mathbf{N}$ be the normal to the surface and $\mathbf{n}$ be the exterior
normal to $\Gamma $ within $S$. Let $B_{\alpha }^{\alpha }$ denote the mean
curvature of $S$ consistent with the direction of $\mathbf{N}$\textbf{\ \cite%
{LeviCivitaTensors}, \cite{McConnellTensors}}. Then%
\begin{equation}
\int_{P}\mathbf{N}B_{\alpha }^{\alpha }dS=\int_{\Gamma }\mathbf{n}d\Gamma .
\label{Integral}
\end{equation}%
In particular, for a patch shrinking to a point, we have%
\begin{equation}
\mathbf{N}B_{\alpha }^{\alpha }=\lim \frac{1}{A}\int_{\Gamma }\mathbf{n}%
d\Gamma .  \label{Limit}
\end{equation}%
The right hand side of equation (\ref{Limit})\ is of coordinate free form
and yields a pure geometric interpretation of mean curvature.

One immediate corollary of equation (\ref{Integral}) is that the integral $%
\int_{\Gamma }\mathbf{n}d\Gamma $ vanishes for any contour in the interior
of a minimal surface.

The derivation is quite brief and we cannot help but wonder whether equation
(\ref{Integral}) could have really been overlooked. The quantity $\mathbf{N}%
B_{\alpha }^{\alpha }$ equals the surface divergence of the surface
covariant basis $\mathbf{S}_{\alpha }$ \cite[equivalent to eqn. (22) on p.
200]{McConnellTensors}:%
\begin{equation}
\mathbf{N}B_{\alpha }^{\alpha }=\nabla ^{\alpha }\mathbf{S}_{\alpha }\mathbf{%
.}
\end{equation}%
By Gauss's theorem on the patch $P$, we have%
\begin{equation}
\int_{P}\nabla ^{\alpha }\mathbf{S}_{\alpha }dS=\int_{\Gamma }n^{\alpha }%
\mathbf{S}_{\alpha }d\Gamma .
\end{equation}%
Recognizing that $n^{\alpha }\mathbf{S}_{\alpha }=\mathbf{n}$, we arrive at
formula (\ref{Integral}).

Formula (\ref{Integral}) suggests a natural way of carrying over the
concepts of the normal vector and mean curvature to triangulated surfaces.
We note in advance that the proposed calculation yields an object analogous
to the \textit{product} $\mathbf{N}B_{\alpha }^{\alpha }$. Therefore, for
minimal surfaces ($\mathbf{N}B_{\alpha }^{\alpha }=\mathbf{0}$) it will not
yield the normal direction. Practical computational problems are almost sure
to arise when computing normals for \textit{nearly} minimal surfaces. On the
other hand, the product $\mathbf{N}B_{\alpha }^{\alpha }$ is a more
fundamental object than either $\mathbf{N}$ or $B_{\alpha }^{\alpha }$, and,
more often than not, $B_{\alpha }^{\alpha }$ arises precisely in the
combination $\mathbf{N}B_{\alpha }^{\alpha }$.

Suppose that on a triangulated surface, triangles $T_{i}$, typically six in
number, with areas $A_{i}$ meet at the node $O$. Let $\mathbf{n}_{i}$ be
unit vectors in the plane of $T_{i}$ perpendicular to the triangle edges
opposite of $O$ and let $\mathbf{n}_{i}$ point away from $O$. Then the
vector quantity%
\begin{equation}
\mathbf{B}=\frac{\sum a_{i}\mathbf{n}_{i}}{\sum A_{i}}  \label{Definition}
\end{equation}%
may be taken as the definition of \textit{vector mean curvature} for
triangulated surfaces. The definition (\ref{Definition}) is easily extended
to the surface Laplacian of general invariant fields by replacing $\mathbf{n}%
_{i}$ with the normal derivative of the field.

The quantity $\mathbf{B}$ has a number of properties in common with its
continuous analogue $\mathbf{N}B_{\alpha }^{\alpha }$, including:

i) $\mathbf{B=0}$ for flat surfaces.

ii) More generally, $\mathbf{B}$ is the gradient of the total area with
respect to the location of the vertices. Therefore, $\mathbf{B=0}$ for
minimal triangulated surfaces.

In conclusion, the central relationships (\ref{Integral}) and (\ref{Limit})
are analogous to Gauss's geometric construction for intrinsic curvature and
may lead to effective computational applications. Surprisingly, we have not
been able to locate these relationships in the existing literature and wish
to share these useful identities as we continue our search.

\bibliographystyle{abbrv}
\bibliography{Classics,ElectronBubbles,FluidFilms,Fluids,InnerCore,Misc,MovingSurfaces,PGrinfeld,Strang,Tensors}

\end{document}